\documentclass[12pt,twoside]{article}
\usepackage{amssymb, amsmath, amsthm,latexsym,graphics,cite}
\usepackage[cp1251]{inputenc}
\usepackage[english]{babel}

\newtheorem{theorem}{Theorem}

\newtheorem{proposition}[theorem]{Proposition}

\newcommand{\seq}[2]{#1_1,#1_2,\ldots ,#1_#2}

\newcommand{\mb}[1]{\mathbb{#1}}
\newcommand{\ov}[1]{\overline{#1}}

\title{On Spectra of Starlike Graphs}
\author{I. K. Redchuk~$^\star$}
\date{}

\begin{document}

\makeatletter
\renewcommand{\@biblabel}[1]{#1.}
\makeatother

\maketitle

\vspace{-15pt}

\noindent
$^\star$~{\small Institute of Mathematics of NAS of Ukraine,\\
$\phantom{^\star}$~Tereschenkovska str., 3, 01601, Kiev, Ukraine\\
$\phantom{^\star}$~E-mail: redchuk@mail.ru}

\bigskip

In the present paper we show that the spectrum of an arbitrary
starlike graph can be completely determined via separating functions
$\rho_t$ (see \cite{NazRoi,RedRoi,Red3}). This fact helps to get in an easy
way several results for the spectra of starlike graphs.

\section{Separating functions and spectra of starlike graphs}

Let $G$ be a simple nonoriented graph, $G_v$ is the set of vertices and
$G_e$ is the set of edges of $G$ ($G_e\subset G_v\times G_v$).
Denote $M(g)=\{h\in G_v\,|\,(h,g)\in G_e\}$. A vertex $g\in G_v$
is called a {\it branching vertex} if $|M(g)|\ge 3$. Graph $G$ is
{\it starlike} if $G$ is a tree and $G$ has at most one branching
vertex. For a starlike graph $G$ the set $G_v$ is
$$
G_v = B_0\sqcup B_1\sqcup B_2 \sqcup \ldots \sqcup B_s,
$$
\noindent where $B_i\cap B_j = \varnothing$, $i,j=\ov{0,s}$,
$B_0=\{g_0\}$, here $g_0$ is a branching vertex in $G$
(if it exists) and $g,h\in B_i$ for some $i$ iff the path of minimal length
beginning in $g$ and ending in $h$ does not contain $g_0$. The sets $B_i$
are {\it branches} of the graph $G$. If $|B_i|=n_i$ then we denote the corresponding
starlike graph by $S_{s;\;\seq{n}{s}}$ or simply $S_{\seq{n}{s}}$.

The {\it spectrum $\sigma(G)$ of a graph} $G$ is a spectrum
of its incidence matrix. Since the last is symmetrical, the spectrum of any graph
is real. Denote $\sigma^+(G)=\{\lambda\in\sigma(G)\,|\,\lambda > 0\}$,
$\sigma_0(G)=\sigma(G)\setminus\{0\}$.
The number $r=\max\{\sigma(G)\}$ is called the index of $G$.
It is known (see, for instance, \cite{CveDoobSac}), that the index
of an arbitrary connected graph $G$ is a simple eigenvalue, and
corresponding eigenspace is spanned by a vector with positive coordinates.
Such vector is called {\it the principal eigenvector of a graph} $G$.

Given $t>0,\, t\in\mb{R}$, let us define a numerical function $\rho_t$:
\begin{equation}\label{eq_rhot}
\rho_t(0)=0,\quad \rho_t(n+1)=t/\left(t-\rho_t(n)\right)
\end{equation}
At that, if $t=\rho_t(n)$ for some $t$ and $n$, then $\rho_t(n+1)$
is not defined; in this case we write $\rho_t(n+1)=\infty$ and set
$\rho_t(n+2)=0$.
Given an integral vector $(\seq{n}{s}),\;n_i\in\mb{N}$
\begin{equation}\label{eq_rhotsum}
\rho_t(\seq{n}{s})=\sum_{i=1}^s\rho_t(n_i).
\end{equation}
Note that in \cite{Red3} functions $\rho_t$ were defined by formulas
\eqref{eq_rhot},\eqref{eq_rhotsum} for $t\ge 1$.

The connection between separating functions $\rho_t$ and indices of starlike graphs
was established in \cite{Red3}. Further we prove more general statement.

Define for $t>0,\; t\in\mb{R}$ the sequence $\{v_n(t)\}$:
\begin{equation}\label{eq_v}
v_0(t)=0,\quad v_1(t)=1,\quad
v_{n+2}(t)=\sqrt{t}v_{n+1}(t)-v_n(t).
\end{equation}
It easy to check (see~\cite{Red3}) that for $n\in\mb{Z}^+$
\begin{equation}\label{eq_vrho}
\rho_t(n) = \sqrt{t}\frac{v_n(t)}{v_{n+1}(t)}.
\end{equation}

\begin{theorem}\label{th_main}
Let $S_{\seq{n}{s}}$ be a starlike graph:
$G_v = B_0\sqcup B_1\sqcup \ldots \sqcup B_s$,
$B_k = \{g_1^k,\ldots,g_{n_k}^k\}$, $k=\ov{1,s}$, $B_0=\{g_0\}$,
and the vertices in branches $B_k$
are numbered in the following way: $g_i^k\in M(g_{i+1}^k)$,
$i=\ov{1,n_k-1}$, $g_{n_k}^k\in M(g_0)$.
\begin{enumerate}
\item For any $t=\lambda^2$, such that $\lambda\in\sigma^+(G)$,
either $t=\rho_t(\seq{n}{s})$ or
there exist such $k,l,\; k\neq l,\; 1\leq k,l\leq s$ that
$t=\rho_t(n_k-1)=\rho_t(n_l-1)$.
\item Conversely: if either $\rho_t(\seq{n}{s})=t$ or there exist
such $k,l,\; k\neq l,\; 1\leq k,l\leq s$ that $\rho_t(n_k-1)=\rho_t(n_l-1)=t$,
then $\sqrt{t}\in\sigma^+(G)$.
\item The index of a graph $G$ equals $\sqrt{t_{\max}}$,
where $t_{\max}$ is the maximal root
of the equation $\rho_t(\seq{n}{s})=t$ in $t$.
\item The components $y_m^k=y(g_m^k),\; y_0=y(g_0)$ of the principal
eigenvector $y$ of a graph $G$ can be calculated
(up to common nonzero multiplier) by the formulas
$$
y_0=\sqrt{t},\quad
y_m^k=t^{(m-n_k)/2}\cdot\prod_{i=m}^{n_k}\rho_t(i), \quad
m=\ov{1,n_k},\; 1\leq k\leq s.
$$
\end{enumerate}
\end{theorem}

\emph{Proof.}

1. Let $\lambda = \sqrt{t}\in\sigma^+(G)$, $x$ is the corresponding
eigenvector, $x_m^k=x(g_m^k),\; x_0=x(g_0)$.
Then the following equalities hold:
\begin{gather}
\lambda x_1^k = x_2^k,\nonumber\\
\lambda x_2^k = x_1^k+x_3^k,\nonumber\\
\lambda x_3^k = x_2^k+x_4^k,\label{eq_eigv}\\
\cdots\; \cdots\; \cdots\; \cdots \nonumber\\
\lambda x_{n_k}^k = x_{n_{k-1}}^k+x_0,\;k=\ov{1,s}\nonumber\\
\lambda x_0 = x_{n_1}^1+\cdots +x_{n_s}^s\nonumber
\end{gather}
Then for all $k=\ov{1,s}$
\begin{gather*}
x_2^k=\lambda x_1,\\
x_{i+2}^k=\lambda x_{i+1}^k-x_i^k,\; i=\ov{1,n_k-2},\\
x_0=\lambda x_{n_k}^k-x_{n_k-1}^k,
\end{gather*}
\noindent hence we obtain that
\begin{gather}
x_{i}^k=v_{i}(t)x_1,\; i=\ov{2,n_k}\nonumber\\
x_0=v_{n_k+1}(t)x_1,\label{eq_vi}
\end{gather}
\noindent here $v_n(t)$ is the sequence defined in \eqref{eq_v}.

If $v_{n_k+1}(t)\neq 0$ for all $k$, then from \eqref{eq_vi} we have
$x_{n_k}^k=\frac{v_{n_k}(t)}{v_{n_k+1}(t)}x_0$.
Then, multiplying the last equality of the system \eqref{eq_eigv}
by $\lambda\neq 0$, we obtain
\begin{equation}\label{eq_alm}
\lambda^2 x_0=\sum_{k=1}^s \frac{\lambda v_{n_k}}{v_{n_k+1}}x_0.
\end{equation}
Equalities \eqref{eq_vi} imply $x_i^k=\frac{v_i(t)}{v_{n_k+1}(t)}x_0$,
$i=\ov{1,n_k},\; k=\ov{1,s}$. Then $x_0\neq 0$, otherwise all components
of eigenvector would vanish. Reducing the both parts of equality \eqref{eq_alm}
by $x_0$ and keeping in mind \eqref{eq_vrho}, we obtain
$$
t=\sum_{k=1}^s\rho_t(n_k)=\rho(\seq{n}{s}).
$$
If $v_{n_k+1}(t)=0$ for some $k$ then \eqref{eq_vrho} implies
$\rho_t(n_k)=\infty$, and so $t=\rho_t(n_k-1)$. Then we obtain $x_0=0$
from \eqref{eq_vi}.
Suppose $\rho_t(n_j-1)\neq t$ for all other $j\neq k,\; 1\leq j\leq s$.
Then $x_i^j=\frac{v_i(t)}{v_{n_j+1}(t)}x_0=0$ for all $j\neq k,\; i=\ov{1,n_j}$,
and due to the last equality from \eqref{eq_eigv} we have $x_{n_k}^k=0$, so
$x_1^k=\frac{x_{n_k}^k}{v_{n_k}(t)}=0$ and, consequently,
$x_i^k=0$ for all $i=\ov{1,n_k}$. As long as all components
of an eigenvector cannot vanish, we come to contradiction; so
there can be found such $l\neq k$, that $\rho_t(n_l-1)=t$.

\medskip

2. Let $t$ be such number that
$\rho_t(n_k-1)=\rho_t(n_l-1)=t,\; k\neq l$ for some $k$, i.~e.
$\rho_t(n_k)=\rho_t(n_l)=\infty$.
Put $x_0=x_i^j=0$ for all such $j=\ov{1,s}$ that $j\neq k$ and $j\neq l$,
$x_i^k=v_{n_k-i+1}(t),\; x_h^l=-v_{n_l-h+1}(t),\;
i=\ov{1,n_k},\; h=\ov{1,n_l}$. It is easy to check that for
$\lambda=\sqrt{t}$ the vector $x$, defined in such a way, satisfies
the equalities \eqref{eq_eigv}, so $\lambda\in\sigma^+(G)$.

Now let $\rho_t(\seq{n}{s})=t$. Then $v_{n_k+1}(t)\neq 0$ for all
$k=\ov{1,s}$. Define in this case the eqigenvector $x$ as follows:
$x_0=1,\; x_i^k=\frac{v_i(t)}{v_{n_k+1}(t)}$. Easily, the equalities
\eqref{eq_eigv} are true again, so $\lambda=\sqrt{t}\in\sigma^+(G)$.

\medskip

3. Let $r$ be the index of a graph $G$, $t=r^2$. Then, by the statement
in item~1, either $\rho_t(\seq{n}{s})=t$ or there exist such $k,l,k\neq l$ that
$\rho_t(n_k-1)=\rho_t(n_l-1)=t$. The second case is impossible since
the component of the eigenvector in the vertex $g_0$ is nonzero
(formulas \eqref{eq_vi}).

\medskip

4. Let $r=\sqrt{t}$  be the index and $y$ be the principal eigenvector
of $G$. Since all its components are nonzero, for any
$k=\ov{1,s}$ we have $v_{n_k+1}\neq 0$, an so, due to \eqref{eq_alm}
and putting $y_0=r$, we obtain
$$
y_m^k=r\frac{v_m(t)}{v_{n_k+1}(t)}=
\frac{1}{r^{n_k-m}}\prod_{i=m}^{n_k}r\frac{v_i(t)}{v_{i+1}(t)}=
\frac{1}{t^{(n_k-m)/2}}\prod_{i=m}^{n_k}\rho_t(i).
$$
This completes the proof of the theorem.

\section{Integral starlike graphs}

A graph $G$ is called {\it integral} if it has integral spectrum. The problem
of description of integral graphs was posed in \cite{HaS2}.
At the present time all integral graphs are obtained only for several
special types of graphs, for example cubic \cite{BuCv}
or with restriction on degrees of vertices of a graph \cite{CvGT3}.
In \cite{CapMau} all integral trees with diameter (i.~e. the longest
of the minimal paths between any two vertices) not exceeding $3$ are described.
For trees with greater diameters it were obtained several (sometimes infinite)
classes, see, for instance, \cite{LiWang,WangLi}.
Further, using theorem~\ref{th_main} we will easily obtain all starlike integral graphs
(it is one of the results of the work \cite{WatSch}).

\begin{proposition}{\rm \cite{RedRoi}}\label{pro_rhoint}
Given $t\in \mb{N},\; t\geq 4$, the list of all integral
vectors $(\seq{n}{s})$ satisfying the equation $\rho_t(\seq{n}{s})=t$
is
\begin{equation}\label{eq_list}
(\underbrace{1,1,\ldots,1}_{t}),\quad
(\underbrace{2,2,\ldots,2}_{t-1}),\quad
(\underbrace{1,3,3,\ldots,3}_{t-1}),\quad
(\underbrace{1,2,5,5,\ldots,5}_{t-1}).
\end{equation}
\end{proposition}

It is known (see \cite{CveDoobSac}) that all graphs with index $r < 2$
(i.~e. Dynkin graphs) are not integral, except graphs $A_1$ (one
vertex) and $A_2$ (two vertices and an edge between them).
Theorem~\ref{th_main} and Proposition~\ref{pro_rhoint} imply that
a starlike graph has integral index $r\geq 2$ if and only if
corresponding integral vector $(\seq{n}{s})$ is contained in the list
\eqref{eq_list} and $t=r^2$.

\begin{proposition}{\rm\cite{CveDoobSac}}\label{pro_sym}
If a graph $G$ is bipartite and $\lambda\in\sigma(G)$ then
$-\lambda\in\sigma(G)$.
\end{proposition}

Using Theorem~\ref{th_main} and Proposition~\ref{pro_sym},
let us check what are remaining elements of spectra of graphs from
the list \eqref{eq_list} for $t\geq 4$.

\medskip

1. Graph $S_{t;\;1,1,\ldots,1}$. The equation $\rho_\tau(\seq{n}{s})=\tau$ in
this case has the form
$t\cdot 1=\tau$, therefore there are no new values besides $\sqrt{t}$. The equation
$\rho_\tau(0)=\tau$ has no solutions for $\tau>0$. Thus, by Theorem~\ref{th_main}
$$
\sigma_0(S_{t;\;1,1,\ldots,1})=\{\sqrt{t},\,-\sqrt{t}\}.
$$

2. Graph $S_{t-1;\;2,2,\ldots,2}$. The equation $\rho_\tau(\seq{n}{s})=\tau$
has the form $(t-1)\cdot \tau/(\tau-1)=\tau$, and
there are no new values besides $\sqrt{t}$ again.
The equation $\rho_\tau(1)=\tau$ has a solution $\tau=1$;
there are more than one branch with $n_k=2$. Consequently,
$$
\sigma_0(S_{t-1;\;2,2,\ldots,2})=\{\sqrt{t},\,-\sqrt{t},\,1,\,-1\}.
$$

3. Graph $S_{t-1;\;1,3,3,\ldots,3}$.
The equation $\rho_\tau(2)=\tau$ has a solution $\tau=2$ and
are more than one branch with $n_k=3$. Therefore,
$$
\sqrt{2}\in \sigma_0(S_{t-1;\;1,3,3,\ldots,3}).
$$

4. Graph $S_{t-1;\;1,2,5,5,\ldots,5}$.
The equation $\rho_\tau(4)=\tau$ has a solution $\tau=3$ and there are more than
one branch with $n_k=5$ when $t>4$. If $t=4$ then the equation
$\rho_\tau(1,2,5)=\tau$ has the form
$1+\frac{\tau}{\tau-1}+\frac{\tau^2-3\tau+1}{(\tau-1)(\tau-3)}=\tau$, which
root is $\tau=(\sqrt{5}+3)/2=\left((\sqrt{5}+1)/2\right)^2$. Therefore,
\begin{eqnarray*}
\sqrt{3}&\in &\sigma_0(S_{t-1;\;1,2,5,5,\ldots,5})\text{ for }t>4,\\
(\sqrt{5}+1)/2&\in& \sigma_0(S_{3;\;1,2,5}).
\end{eqnarray*}

Thus, we proved the following

\begin{theorem}\label{pro_intstarshap}{\rm \cite{WatSch}}.
All integral starlike graphs are $A_1$,
$S_{t;\;1,1,\ldots,1}$ and
$S_{t-1;\;2,2,\ldots,2}$,
where $t=r^2$ for some $r\in \mb{N}$.
\end{theorem}

\end{document}